# LIFTING ENHANCED FACTORIZATION SYSTEMS TO FUNCTOR 2-CATEGORIES

PETER J. HAINE

ABSTRACT. In this paper we give sufficient conditions for lifting an enhanced factorization system $(\mathcal{E}, \mathcal{M})$ on a 2-category $D$ to the functor 2-category $D^C$, where $C$ is a small 2-category. Due to previous work of Lack, our work provides coherence results for 2-monads on functor 2-categories. These coherence results are of immediate interest to an ongoing collaboration of Guillou, May, Merling, and Osorno on *equivariant infinite loop space theory*.

## 0. Overview

The motivation for this work comes from coherence theory for pseudo-algebras over 2-monads, specifically 2-monads on functor 2-categories. Our work is of immediate interest to an ongoing collaboration of Bert Guillou, Peter May, Mona Merling, and Angélica Osorno on *equivariant infinite loop space theory* [3]. The basic problem is to provide a setting where we can answer the following question: given a 2-category $D$ and a 2-monad $T$ on $D$, is every pseudo-$T$-algebra equivalent (as a pseudo-$T$-algebra) to a strict $T$-algebra? In other words, this question asks if every pseudo-$T$-algebra be "strictified" to a strict $T$-algebra. Given any small set $S$ (regarded as a discrete 2-category) and a 2-monad on the functor 2-category $\mathbf{Cat}^S$, Power [7, Cor. 3.5] showed that if $T$ preserves 1-cells in $\mathbf{Cat}^S$ that are $S$-indexed sets of bijective-on-objects functors, then every pseudo-$T$-algebra is equivalent to a strict one. There are two ideas that Power employed in his proof; the first is to lift the result "levelwise" to the functor 2-category, and the second is that the preservation of a certain class of 1-cells is a sufficient condition for "strictification".

Lack [6, §4.2] used the notion of an *enhanced factorization system,* originally developed by Kelly [4], to generalize (and strengthen) Power's result to a 2-category $A$ with an enhanced factorization system satisfying one additional condition, which we call *rigidity*. The motivation for our work is to generalize Power's result to the functor 2-category $\mathbf{Cat}^C$, or, more generally, to the functor 2-category $D^C$, where $D$ is a 2-category and $C$ is a small 2-category. In light of Lack's work [6, §4.2], this generalization can be achieved by specifying a rigid enhanced factorization system on $D^C$. Without an enhanced factorization system on $D$, specifying an enhanced factorization system on $D^C$ is hopeless, hence the goal of this paper is to provide sufficient conditions to lift an enhanced factorization system on a 2-category $D$ to an enhanced factorization system on the functor 2-category $D^C$, defined "levelwise".

In §1 we review enhanced factorization systems and precisely state Lack's generalization of Power's result. We also show how to lift the factorization of 1-cells in a 2-category $D$ with an enhanced factorization system to a levelwise factorization of 1-cells in the functor 2-category $D^C$, where $C$ is a small 2-category as this requires no additional assumptions on the enhanced factorization system. In §2 we analyze the enhanced factorization system on $\mathbf{Cat}$ in order to determine sufficient conditions for lifting an enhanced factorization system on a 2-category $D$ to $D^C$. We conclude §2 by stating the main result of this paper, which says that, under mild assumptions, an enhanced factorization system on a 2-category $D$ can be lifted to an enhanced factorization system on $D^C$, in particular, our results hold when $D = \mathbf{Cat}$. Sections 3 and 4 are dedicated to proving this result.

**Acknowledgments.** The author completed much of this work while at the University of Chicago REU program funded by NSF DMS-1344997. He would like to thank Peter May for asking him to work on this project, for much of the terminology used here, his advice, and many helpful comments. He would also like to thank Emily Riehl for teaching him most of the category theory he knows, first suggesting that he read about 2-categories, and answering his many questions

## 1. Enhanced Factorization Systems

We begin by reviewing *enhanced factorization systems.* Before we do so, let us first set some notational conventions.

1.1. **Notation.** Since we need to go back and forth between general definitions and specific instances of the definitions, where there is a preferred classical choice of notation, below we outline our general notational conventions for the convenience of the reader.

— We write **Cat** for the 2-category of small categories, functors, and natural transformations.
— Throughout, the boldface Roman letters $A$, $C$, and $D$ denote 2-categories.
— We always use the letter $C$ for a small 2-category. We denote the objects of $C$ by lowercase Roman letters $c$, $c'$, etc., 1-cells by lowercase Roman letters such as $f$, $g$, or $h$, and 2-cells by lowercase Greek letters such as $\lambda$ and $\eta$.
— We always use the letter $A$ for an arbitrary 2-category: since we are mostly interested in the case that $A$ is the functor 2-category $D^C$, we write capital Roman letters $F$, $G$, etc. for objects of a 2-category, lowercase Greek letters $\alpha$, $\beta$, etc. for 1-cells, and uppercase Greek letters $\Psi$, $\Phi$, etc. for 2-cells of $A$.
— We use calligraphic letters such as $\mathcal{E}$ and $\mathcal{M}$ to denote distinguished classes of 1-cells in a 2-category.
— We use the double-struck arrow "$\Longrightarrow$" to denote a 2-cell in a 2-category, for example, a natural transformation in **Cat**.

1.2. **Definition** ([6, § 4.2]). An *enhanced factorization system* on a 2-category $A$ consists of a pair of classes of 1-cells $(\mathcal{E}, \mathcal{M})$, both containing all isomorphisms, satisfying the following properties.[1]

(1.2.a) Every 1-cell $\alpha$ of $A$ factors (not necessarily uniquely) as a composite $\alpha = \mu \circ \varepsilon$, where $\varepsilon \in \mathcal{E}$ and $\mu \in \mathcal{M}$.

(1.2.b) For a diagram in $A$ of the form

$$\begin{array}{ccc} F & \xrightarrow{\varepsilon} & F' \\ \alpha \downarrow & \Swarrow_{\Psi} & \downarrow \alpha' \\ G & \xrightarrow[\mu]{} & G' \end{array}$$

---

[1]The notations "$\mathcal{E}$" and "$\mathcal{M}$" are meant to stand for "epimorphism" and "monomorphism", respectively. Example 1.4 below shows that 1-cells in $\mathcal{E}$ are not necessarily epimorphisms, and 1-cells in $\mathcal{M}$ are not necessarily monomorphism, but in §2 we will see that in cases of interest they will satisfy properties with respect to 2-cells reminiscent of epimorphisms and monomorpisms, respectively.

where $\varepsilon \in \mathcal{E}$, $\mu \in \mathcal{M}$, and $\Psi$ is an invertible 2-cell, there is a unique pair $(\delta, \widetilde{\Psi})$ consisting of a 1-cell $\delta \colon F' \longrightarrow G$ and an invertible 2-cell $\widetilde{\Psi} \colon \alpha' \Longrightarrow \mu\delta$ so that we have a factorization

$$\begin{array}{ccc} F \xrightarrow{\varepsilon} F' & & F \xrightarrow{\varepsilon} F' \\ \alpha \downarrow \; \Swarrow_{\Psi} \; \downarrow \alpha' & = & \alpha \downarrow \; {}^{\delta}\!\!\Swarrow_{\overline{\Psi}} \; \downarrow \alpha' \\ G \xrightarrow{\mu} G' & & G \xrightarrow{\mu} G' \, , \end{array}$$

where $\delta\varepsilon = \alpha$. By the uniqueness of $\widetilde{\Psi}$, if $\Psi$ is the identity, then $\gamma\delta = \alpha'$ and $\widetilde{\Psi}$ is the identity.

(1.2.c) Suppose that we are given 1-cells $\varepsilon \colon F \longrightarrow F'$ in $\mathcal{E}$ and $\mu \colon G \longrightarrow G'$ in $\mathcal{M}$, two parallel pairs $\alpha_1, \alpha_2 \colon F \Longrightarrow G$ and $\alpha'_1, \alpha'_2 \colon F' \Longrightarrow G'$ in $A$, so that the squares

$$\begin{array}{ccc} F \xrightarrow{\varepsilon} F' & & F \xrightarrow{\varepsilon} F' \\ \alpha_1 \downarrow \qquad \downarrow \alpha'_1 & \text{and} & \alpha_2 \downarrow \qquad \downarrow \alpha'_2 \\ G \xrightarrow{\mu} G' & & G \xrightarrow{\mu} G' \end{array}$$

commute, and 2-cells $\Phi \colon \alpha_1 \Longrightarrow \alpha_2$ and $\Phi' \colon \alpha'_1 \Longrightarrow \alpha'_2$ such that $\mu\Phi = \Phi'\varepsilon$. Let $\delta_1$ and $\delta_2$ denote the unique 1-cells given by (1.2.b) so that each of the sub-triangles in the diagrams

$$\begin{array}{ccc} F \xrightarrow{\varepsilon} F' & & F \xrightarrow{\varepsilon} F' \\ \alpha_1 \downarrow \; {}^{\delta_1}\!\!\nearrow \; \downarrow \alpha'_1 & \text{and} & \alpha_2 \downarrow \; {}^{\delta_2}\!\!\nearrow \; \downarrow \alpha'_2 \\ G \xrightarrow{\mu} G' & & G \xrightarrow{\mu} G' \end{array}$$

commute. Then there exists a unique 2-cell $\Delta \colon \delta_1 \Longrightarrow \delta_2$ so that both $\Delta\varepsilon = \Phi$ and $\mu\Delta = \Phi'$.

We say that an enhanced factorization system $(\mathcal{E}, \mathcal{M})$ is *rigid* if the following additional property holds.

(1.2.d) Given 1-cells $\mu \colon F \rightleftarrows G \colon \alpha$ in $A$, where $\mu \in \mathcal{M}$, if $\mu\alpha \cong \mathrm{id}_G$, then $\alpha\mu \cong \mathrm{id}_F$.

1.3. **Remark.** In [6, §4.2] when Lack defines an enhanced factorization system on a 2-category $D$, he requires that the pair $(\mathcal{E}, \mathcal{M})$ be a "factorization system" on the underlying 1-category $D_1$ of $D$. However, there does not seem to be consistency in the literature about the meaning of a "factorization system". There is an ambiguity as to whether a "factorization system" should refer to an "orthogonal factorization system", a "weak factorization system", the condition (1.2.b) of Definition 1.2, or something else entirely. Since Lack [6] does not specify which of notions he means, we have chosen to adopt the least restrictive notion.

1.4. **Example.** The following classical examples of enhanced factorization systems are essentially due to Power [7].

(1.4.a) The 2-category **Cat** has an rigid enhanced factorization system $(\mathcal{B}, \mathcal{F})$ where $\mathcal{B}$ is the class of bijective-on-objects functors, and $\mathcal{F}$ is the class of fully faithful functors.

(1.4.b) If $S$ is a small set, the 2-category $\mathbf{Cat}^S$ has an rigid enhanced factorization system $(\mathcal{B}^S, \mathcal{F}^S)$ where $\mathcal{B}^S$ is the class of $S$-indexed sets of bijective-on-objects functors and $\mathcal{F}^S$ is the class of $S$-indexed sets of fully faithful functors.

(1.4.c) More generally, if $S$ is a small set, and $\boldsymbol{D}$ is a 2-category with a rigid enhanced factorization system $(\mathcal{E}, \mathcal{M})$, then $\boldsymbol{D}^S$ has an rigid enhanced factorization system $(\mathcal{E}^S, \mathcal{M}^S)$ where $\mathcal{E}^S$ is the class of $S$-indexed sets of 1-cells in $\mathcal{E}$, and $\mathcal{M}^S$ is the class of $S$-indexed sets of 1-cells in $\mathcal{M}$.

Coherence results for 2-monads on each of the 2-categories described in Example 1.4 are immediate consequences of Lack's general coherence result [6, §4.2], which we explain presently. We assume familiarity with the theory of 2-monads and do not review these notions as the main results and proofs in this paper do not actually require any knowledge of 2-monads (although 2-monads do form the motivation for our work). The unfamiliar reader should consult [1; 5, §§3.1–3.2; 6, §1; 7, §2] for overviews of the basic theory.

1.5. **Notation.** Suppose that $A$ is a 2-category and that $T$ is a 2-monad on $A$. Write $\mathbf{Alg}_T$ for the 2-category of (strict) $T$-algebras, morphisms, and algebra 2-cells. Write $\mathbf{Alg}_T^{ps}$ for the 2-category of pseudo-$T$-algebras, morphisms, and algebra 2-cells.

1.6. **Remark.** In [1, §1.2] Blackwell, Kelly, and Power discuss the different notions of strict, lax, and pseudo-$T$-algebras: in their notation, the 2-category $\mathbf{Alg}_T^{ps}$ is written as $T$-Alg. Power [7, §2] also has a detailed discussion of pesudo-$T$-algebras: in his paper he writes $\mathbf{PS\text{-}T\text{-}Alg}$ for the 2-category $\mathbf{Alg}_T^{ps}$.

1.7. **Definition.** Suppose that $A$ is a 2-category and that $\mathcal{E}$ is a class of 1-cells in $A$. We say that a 2-monad $T$ on $A$ *preserves* $\mathcal{E}$ if $\varepsilon \in \mathcal{E}$ implies that $T\varepsilon \in \mathcal{E}$.

Lack's generalization of Power's result is the following.

1.8. **Theorem** ([6, Thm. 4.6]). *Suppose that $A$ is a 2-category with a rigid enhanced factorization system $(\mathcal{E}, \mathcal{M})$, and that $T$ is a 2-monad on $A$. If $T$ preserves $\mathcal{E}$, then the inclusion $\mathbf{Alg}_T \hookrightarrow \mathbf{Alg}_T^{ps}$ has a left 2-adjoint, and the components of the unit of the adjunction are equivalences in $\mathbf{Alg}_T^{ps}$.*

In light of Lack's result, one way of generating coherence results for 2-monads on the functor 2-category $\boldsymbol{D}^{\boldsymbol{C}}$, where $\boldsymbol{D}$ is a 2-category with an enhanced factorization system $(\mathcal{E}, \mathcal{M})$ and $\boldsymbol{C}$ is a small 2-category, is to show that $\boldsymbol{D}^{\boldsymbol{C}}$ has a rigid enhanced factorization system. The first step toward tackling this problem is to simply show any 1-cell in $\boldsymbol{D}^{\boldsymbol{C}}$ factors as a composite of a 1-cell in a specified class $\mathcal{E}^{\boldsymbol{C}}$, followed by a 1-cell in another specified class $\mathcal{M}^{\boldsymbol{C}}$; this is what we are concerned with for the rest of this section. We can prove this result with our bare hands, using no additional assumptions, however, we shall soon see that the other axioms require more care.

1.9. **Proposition.** *Suppose that $\boldsymbol{C}$ is a small 2-category and that $\boldsymbol{D}$ is a 2-category with an enhanced factorization system $(\mathcal{E}, \mathcal{M})$. Given objects $F, G \in \boldsymbol{D}^{\boldsymbol{C}}$ and a 1-cell $\alpha \colon F \longrightarrow G$ in $\boldsymbol{D}^{\boldsymbol{C}}$, the 1-cell $\alpha$ factors as $\alpha = \mu \circ \varepsilon$, where, for each $c \in C$, the component $\varepsilon_c$ is in $\mathcal{E}$ and the component $\mu_c$ is in $\mathcal{M}$.*

*Proof.* The goal is to factor the 2-natural transformation $\alpha$ as a composite of 2-natural transformations

$$\begin{array}{c} F \xrightarrow{\alpha} G \\ {}_{\varepsilon} \searrow \quad \nearrow {}_{\mu} \\ I \end{array},$$

through a 2-functor $I \colon \boldsymbol{C} \longrightarrow \boldsymbol{D}$, where the components of $\varepsilon$ lie in $\mathcal{E}$ and the components of $\mu$ lie in $\mathcal{M}$. Since the factorization of 1-cells in $\boldsymbol{D}$ as a composite of a 1-cell in $\mathcal{E}$ followed by a 1-cell in $\mathcal{M}$ is not unique, the factorization through a 2-functor $I$ will not be uniquely defined; instead, we use the properties of the enhanced factroization system on $D$ to choose such a 2-functor $I$.

To define a 2-functor $I\colon C \longrightarrow D$ on objects, for each $c \in C$, use the enhanced factorization system on $D$ to choose a factorization

$$\begin{array}{ccc} F(c) & \xrightarrow{\alpha_c} & G(c) \\ & \searrow^{\varepsilon_c} \quad \nearrow^{\mu_c} & \\ & I(c) & \end{array},$$

where $I(c)$ is an object of $D$, $\varepsilon_c \in \mathcal{E}$, and $\mu_c \in \mathcal{M}$, and define $I$ on objects of $C$ by the assignment $c \longmapsto I(c)$. To define $I$ on 1-cells, we use condition (1.2.b) of Definition 1.2: because $(\mathcal{E}, \mathcal{M})$ is an enhanced factorization system on $D$, for each morphism $f\colon c \longrightarrow c'$ in $C$, there is a unique factorization

$$\begin{array}{ccc} \begin{array}{ccc} F(c) & \xrightarrow{\varepsilon_c} & I(c) \\ F(f)\downarrow & & \downarrow\mu_c \\ F(c') & & G(c) \\ \varepsilon_{c'}\downarrow & & \downarrow G(f) \\ I(c') & \xrightarrow[\mu_{c'}]{} & G(c') \end{array} & = & \begin{array}{ccc} F(c) & \xrightarrow{\varepsilon_c} & I(c) \\ F(f)\downarrow & \searrow & \downarrow\mu_c \\ F(c') & I(f) & G(c) \\ \varepsilon_{c'}\downarrow & \searrow & \downarrow G(f) \\ I(c') & \xrightarrow[\mu_{c'}]{} & G(c') \end{array} \end{array} \tag{1.9.1}$$

so that each of the sub-triangles in the right-hand diagram of (1.9.1) commutes. With this suggestively chosen notation, we define $I$ on 1-cells of $C$ by sending $f$ to the (unique) 1-cell $I(f)$.

Before defining $I$ on the 2-cells of $C$, let us first show the assignment

$$[f\colon c \longrightarrow c'] \longmapsto [I(f)\colon I(c) \longrightarrow I(c')]$$

satisfies the 1-categorical properties of a 2-functor, i.e., defines a functor from the underlying 1-category of $C$ to the underlying 1-category of $D$. To see this, first notice that by the uniqueness of the factorization (1.9.1), it is clear that $I(\mathrm{id}_c) = \mathrm{id}_{I(c)}$. Then by the commutativity of the triangles in the right-hand diagram of (1.9.1), we get a diagram

$$\begin{array}{ccccc} F(c) & \xrightarrow{F(f)} & F(c') & \xrightarrow{F(g)} & F(c'') \\ \varepsilon_c\downarrow & & \downarrow\varepsilon_{c'} & & \downarrow\varepsilon_{c''} \\ I(c) & \xrightarrow{I(f)} & I(c') & \xrightarrow{I(g)} & I(c'') \\ \mu_c\downarrow & & \downarrow\mu_{c'} & & \downarrow\mu_{c''} \\ G(c) & \xrightarrow[G(f)]{} & G(c') & \xrightarrow[G(g)]{} & G(c'') \end{array}, \tag{1.9.2}$$

in which each of the sub-squares commutes. The commutativity of each of the sub-squares of (1.9.2), along with the functoriality of $F$ and $G$ show that the diagram

$$\begin{array}{ccc} F(c) & \xrightarrow{F(gf)} & F(c'') \\ \varepsilon_c \downarrow & & \downarrow \varepsilon_{c''} \\ I(c) & \xrightarrow{I(g)I(f)} & I(c'') \\ \mu_c \downarrow & & \downarrow \mu_{c''} \\ G(c) & \xrightarrow[G(gf)]{} & G(c'') \end{array}$$

commutes, so, by the uniqueness of $I(gf)$, we see that $I(gf) = I(g)I(f)$. Then, by construction, the 1-cells $(\varepsilon_c)_{c \in C}$ and $(\mu_c)_{c \in C}$ satisfy the following 1-categorical naturality conditions: for any 1-cell $f \colon c \longrightarrow c'$ of $C$, the squares

$$\begin{array}{ccc} F(c) & \xrightarrow{F(f)} & F(c') \\ \varepsilon_c \downarrow & & \downarrow \varepsilon_{c'} \\ I(c) & \xrightarrow[I(f)]{} & I(c') \end{array} \quad \text{and} \quad \begin{array}{ccc} I(c) & \xrightarrow{I(f)} & I(c') \\ \mu_c \downarrow & & \downarrow \mu_{c'} \\ G(c) & \xrightarrow[G(f)]{} & G(c') \end{array} \qquad (1.9.3)$$

commute.

Finally let us define $I \colon C \longrightarrow D$ on 2-cells; to do this, we use condition (1.2.c) of Definition 1.2. Suppose that we are given the data of objects, 1-cells, and 2-cells in $C$, as displayed below:

$$c \xrightarrow[g]{\overset{f}{\Downarrow \lambda}} c' \ .$$

In this case, we know that the diagrams

$$\begin{array}{ccc} F(c) & \xrightarrow{\varepsilon_c} & I(c) \\ F(f) \downarrow & & \downarrow \mu_c \\ F(c') & & G(c) \\ \varepsilon_{c'} \downarrow & & \downarrow G(f) \\ I(c') & \xrightarrow[\mu_{c'}]{} & G(c') \end{array} \quad \text{and} \quad \begin{array}{ccc} F(c) & \xrightarrow{\varepsilon_c} & I(c) \\ F(g) \downarrow & & \downarrow \mu_c \\ F(c') & & G(c) \\ \varepsilon_{c'} \downarrow & & \downarrow G(g) \\ I(c') & \xrightarrow[\mu_{c'}]{} & G(c') \end{array}$$

both commute. We also have a pair of (whiskered) two cells

$$\varepsilon_{c'} F(\lambda) \colon \varepsilon_{c'} \circ F(f) \Longrightarrow \varepsilon_{c'} \circ F(g) \quad \text{and} \quad G(\lambda)\mu_c \colon G(f) \circ \mu_c \Longrightarrow G(g) \circ \mu_c \ .$$

Since $\alpha_c$ and $\alpha_{c'}$ factor as $\alpha_c = \mu_c \varepsilon_c$ and $\alpha_{c'} = \mu_{c'} \varepsilon_{c'}$ we know that

$$F(c) \underset{F(g)}{\overset{F(f)}{\Rightarrow}} F(\lambda)\Downarrow F(c') \xrightarrow{\varepsilon_{c'}} I(c') \xrightarrow{\mu_{c'}} G(c') \quad = \quad F(c) \underset{F(g)}{\overset{F(f)}{\Rightarrow}} F(\lambda)\Downarrow F(c') \xrightarrow{\alpha_{c'}} G(c')$$

and

$$F(c) \xrightarrow{\varepsilon_c} I(c) \xrightarrow{\mu_c} G(c) \underset{G(g)}{\overset{G(f)}{\Rightarrow}} \Downarrow G(\lambda)\, G(c') \quad = \quad F(c) \xrightarrow{\alpha_c} G(c) \underset{G(g)}{\overset{G(f)}{\Rightarrow}} \Downarrow G(\lambda)\, G(c').$$

Because $\alpha$ is a 2-natural transformation we also know that

$$F(c) \underset{F(g)}{\overset{F(f)}{\Rightarrow}} \Downarrow F(\lambda)\, F(c') \xrightarrow{\alpha_{c'}} G(c') \quad = \quad F(c) \xrightarrow{\alpha_c} G(c) \underset{G(g)}{\overset{G(f)}{\Rightarrow}} \Downarrow G(\lambda)\, G(c'),$$

which verifies the conditions necessary to apply (1.2.c) of Definition 1.2, which (by our definitions of $I(f)$ and $I(g)$) gives a unique 2-cell $I(\lambda)\colon I(f) \implies I(g)$ in $D$ so that

$$F(c) \underset{F(g)}{\overset{F(f)}{\Rightarrow}} \Downarrow F(\lambda)\, F(c') \xrightarrow{\varepsilon_{c'}} I(c') \quad = \quad F(c) \xrightarrow{\varepsilon_c} I(c) \underset{I(g)}{\overset{I(f)}{\Rightarrow}} \Downarrow I(\lambda)\, I(c') \qquad (1.9.4)$$

and

$$I(c) \underset{I(g)}{\overset{I(f)}{\Rightarrow}} \Downarrow I(\lambda)\, I(c') \xrightarrow{\mu_{c'}} G(c') \quad = \quad I(c) \xrightarrow{\mu_c} G(c) \underset{G(g)}{\overset{G(f)}{\Rightarrow}} \Downarrow G(\lambda)\, G(c'), \qquad (1.9.5)$$

What the equalities of pasting diagrams in (1.9.4) and (1.9.5) tell us is that once we verify that the assignment

$$c \underset{g}{\overset{f}{\Rightarrow}} \Downarrow \lambda\, c' \quad \longmapsto \quad I(c) \underset{I(g)}{\overset{I(f)}{\Rightarrow}} \Downarrow I(\lambda)\, I(c') \qquad (1.9.6)$$

actually defines a 2-functor, the collections of 1-cells $(\varepsilon_c)_{c\in C}$ and $(\mu_c)_{c\in C}$ will satisfy the remaining 2-categorical conditions necessary in order to define a 2-natural transformation.

Therefore, all that remains to be shown is that the assignment of $I$ on 2-cells respects identities and composition. The fact that $I(\mathrm{id}_f) = \mathrm{id}_{I(f)}$ for any 1-cell $f$ of $C$ follows immediately from the uniqueness property of $I(\mathrm{id}_f)$ along with the facts that $F(\mathrm{id}_f) = \mathrm{id}_{F(f)}$ and $G(\mathrm{id}_f) = \mathrm{id}_{G(f)}$. Similarly, given a composable pair of 2-cells in $C$ as displayed below

$$c \underset{h}{\overset{f}{-g-}} \Downarrow \lambda \atop \Downarrow \eta\, c',$$

by the 2-functoriality of $F$ and $G$ we know that $F(\eta)F(\lambda) = F(\eta\lambda)$ and $G(\eta)G(\lambda) = G(\eta\lambda)$, which tells us that there are equalities of whiskered 2-cells

$$\varepsilon_{c'}(F(\eta)F(\lambda)) = \varepsilon_{c'}F(\eta\lambda) \quad \text{and} \quad (G(\eta)G(\lambda))\mu_c = G(\eta\lambda)\mu_c. \qquad (1.9.7)$$

The equalities displayed in (1.9.7) imply that the 2-cell $I(\eta)I(\lambda)\colon I(f) \implies I(h)$ satisfies the uniqueness property of defining $I(\eta\lambda)$, hence $I(\eta)I(\lambda) = I(\eta\lambda)$, as desired. This completes the lengthy verification that the assignments (1.9.6) define a 2-functor $I\colon C \longrightarrow D$. Hence, every 1-cell $\alpha$ in $\boldsymbol{D}^{\boldsymbol{C}}$ factors as a composite $\alpha = \mu\varepsilon$, where all of the components of $\varepsilon$ are in $\mathcal{E}$, and all of the components of $\mu$ are in $\mathcal{M}$. □

The point of this is that an enhanced factorization system on $\boldsymbol{D}$ yields a way of factoring 1-cells in $\boldsymbol{D}^{\boldsymbol{C}}$, defined "levelwise".

1.10. **Definition.** Suppose that $\boldsymbol{D}$ is a 2-category with an enhanced factorization system $(\mathcal{E}, \mathcal{M})$, and that $\boldsymbol{C}$ is a small 2-category. Let $\mathcal{E}^{\boldsymbol{C}}$ denote the collection of 1-cells in $\boldsymbol{D}^{\boldsymbol{C}}$, all of whose components are in $\mathcal{E}$, and $\mathcal{M}^{\boldsymbol{C}}$ denote the collection of 1-cells in $\boldsymbol{D}^{\boldsymbol{C}}$, all of whose components are in $\mathcal{M}$.

Given a 1-cell $\alpha \in \boldsymbol{D}^{\boldsymbol{C}}$, we call a factorization $\alpha = \mu \circ \varepsilon$, where $\varepsilon \in \mathcal{E}^{\boldsymbol{C}}$ and $\mu \in \mathcal{M}^{\boldsymbol{C}}$ (described in Proposition 1.9) a *levelwise factorization* of $\alpha$.

The rest of this paper is concerned with giving sufficient conditions on an enhanced factorization system $(\mathcal{E}, \mathcal{M})$ on a 2-category $\boldsymbol{D}$ so that for any small 2-category $\boldsymbol{C}$, the pair of classes of 1-cells $(\mathcal{E}^{\boldsymbol{C}}, \mathcal{M}^{\boldsymbol{C}})$ defines an enhanced factorization system on $\boldsymbol{D}^{\boldsymbol{C}}$.

## 2. The Enhanced Factorization System on **Cat**

In this section we analyze the enhanced factorization system $(\mathcal{B}, \mathcal{F})$ on **Cat** of bijective-on-objects and fully faithful functors. Many of the results in this section are known, but we provide proofs of them as they motivate generalizations of these notions to arbitrary 2-categories. In addition to motivating the definitions presented in this section, these results also serve to provide an example of our main result, Theorem 2.9, in the case that $\boldsymbol{D} = $ **Cat**. Moreover, they provide intuition about how to generalize the relevant results for **Cat** to the 2-category of categories internal to a fixed bicomplete cartesian monoidal category, which is the primary case of interest to Guillou, May, Merling, and Osorno. For the following few results we use more classical 1-categorical notation.

2.1. **Notation.** In Lemmas 2.2 and 2.5, we write:
   — $C$, $D$, and $E$ for ordinary 1-categories,
   — $B$ for a functor which is bijective-on-objects,
   — $F$ for a fully faithful functor,
   — $G$ and $H$ for arbitrary functors,
   — and lowercase Greek letters such as $\lambda$, $\lambda'$, and $\eta$ for natural transformations.

2.2. **Lemma.** *Suppose that we have categories, functors, and natural transformations as displayed below*

$$C \underset{H}{\overset{G}{\rightrightarrows}} {\lambda\Downarrow \Downarrow\lambda'} D \xrightarrow{F} E,$$

*where $F$ is fully faithful.*

(2.2.a) *If $F\lambda = F\lambda'$, then $\lambda = \lambda'$.*

(2.2.b) *A natural transformation $\eta\colon FG \implies FH$ lifts to a unique natural transformation $\hat{\eta}\colon G \implies H$ with the property that $F\hat{\eta} = \eta$. Moreover, if $\eta$ is a natural isomorphism, so is $\hat{\eta}$.*

(2.2.c) *If $G(c) = H(c)$ for all $c \in C$ and $FG = FH$, then $G = H$. Thus, if $B\colon C' \longrightarrow C$ is a bijective-on-objects functor so that $GB = BH$ and $FG = FH$, then $G = H$.*

*Proof.* First, (2.2.a) follows immediately from the fact that $F$ is fully faithful and for all $c \in C$ the components $\lambda_c$ and $\lambda'_c$ have the same source and target.

Second, suppose that we have a natural transformation $\eta\colon FG \Longrightarrow FH$. Since $F$ is fully faithful, for each $c \in C$, there exists a unique morphism $\hat{\eta}_c\colon G(c) \longrightarrow H(c)$ in $D$ so that $F(\hat{\eta}_c) = \eta_c$. To see that the morphisms $(\hat{\eta}_c)_{c \in C}$ assemble into a natural transformation $\hat{\eta}\colon G \Longrightarrow H$, suppose that $f\colon c \longrightarrow c'$ is a morphism of $C$, and consider the square

$$\begin{array}{ccc} G(c) & \xrightarrow{G(f)} & G(c') \\ \hat{\eta}_c \downarrow & & \downarrow \hat{\eta}_{c'} \\ H(c) & \xrightarrow{H(f)} & H(c') \,. \end{array} \qquad (2.2.1)$$

To see that (2.2.1) commutes, notice that by the functoriality of $F$ and naturality of $\eta$ we have

$$F(\hat{\eta}_{c'} \circ G(f)) = \eta_{c'} \circ FG(f) = FH(f) \circ \eta_c = F(H(f) \circ \hat{\eta}_c) \,.$$

Since $F$ is fully faithful and $\hat{\eta}_{c'} \circ G(f)$ and $H(f) \circ \hat{\eta}_c$ have the same source and target, this implies that $\hat{\eta}_{c'} \circ G(f) = H(f) \circ \hat{\eta}_c$, as desired. Moreover, $\hat{\eta}$ is defined so that $F\hat{\eta} = \eta$. Lastly, since fully faithful functors reflect isomorphisms, if $\eta$ is a natural isomorphism, then so is $\hat{\eta}$.

To prove (2.2.c), all that needs to be verified is that $G(f) = H(f)$ for all morphisms $f$ of $C$. Since $G(c) = H(c)$ for all $c \in C$, the morphisms $G(f)$ and $H(f)$ have the same source and target. Since $FG(f) = FH(f)$, and $F$ is fully faithful, this implies that $G(f) = H(f)$. □

The following definition is a generalization of the situation occurring in item (2.2.b) of Lemma 2.2.

**2.3. Definition.** Suppose that $A$ is a 2-category and that $\mu$ is a 1-cell in $A$. We say that *post-composition with $\mu$ creates invertible 2-cells* if whenever we have 1-cells $\alpha, \beta \in A$ with the same source and target equipped with an invertible 2-cell $\Psi\colon \mu\alpha \Longrightarrow \mu\beta$, there exists a unique invertible 2-cell $\widehat{\Psi}\colon \alpha \Longrightarrow \beta$ so that $\mu\widehat{\Psi} = \Psi$.

Suppose that $\mathcal{M}$ is a class of 1-cells in a 2-category $A$. We say that *post-composition with 1-cells in $\mathcal{M}$ creates invertible 2-cells*, if post-composition with every 1-cell $\mu \in \mathcal{M}$ creates invertible 2-cells.

The following definition is a generalization of the situation occurring in item (2.2.c) of Lemma 2.2.

**2.4. Definition.** We say that an enhanced factorization system $(\mathcal{E}, \mathcal{M})$ on a 2-category $A$ *separates parallel pairs* if whenever we are given a parallel pair $\alpha, \beta\colon F \rightrightarrows G$ so that there exists a 1-cell $\varepsilon \in \mathcal{E}$ with target $F$ such that $\alpha\varepsilon = \beta\varepsilon$ and a 1-cell $\mu \in \mathcal{M}$ such that $\mu\alpha = \mu\beta$, we have $\alpha = \beta$.

Now for a result regarding functors which are bijective-on-objects.

**2.5. Lemma.** *Suppose that we have categories, functors, and natural transformations as displayed below*

$$C \xrightarrow{B} D \underset{H}{\overset{G}{\rightrightarrows}} E \,, \quad \lambda \Downarrow \Downarrow \lambda'$$

*where $B$ is bijective-on-objects. If $\lambda B = \lambda' B$, then $\lambda = \lambda'$.*

*Proof.* Since $\lambda B = \lambda' B$, for all $c \in C$ we have $\lambda_{B(c)} = \lambda'_{B(c)}$. Since $B$ is bijective-on-objects, this says that $\lambda_d = \lambda'_d$ for each $d \in D$, so $\lambda = \lambda'$. □

In a 2-category, 1-cells which are 2-*epimorphisms* express the conclusion of Lemma 2.5 about functors which are bijective-on-objects. Similarly, 2-*monomorphisms* express the dual property of fully faithful functors, which is (2.2.a) of Lemma 2.2.

2.6. **Definition.** We say that a 1-cell $\varepsilon \colon F \longrightarrow G$ in a 2-category $A$ is a 2-*epimorphism* if whenever we have objects, 1-cells, and 2-cells as displayed below

$$F \xrightarrow{\varepsilon} G \overset{\Psi \Downarrow \, \Downarrow \Psi'}{\rightrightarrows} H \,,$$

if $\Psi \varepsilon = \Psi' \varepsilon$, then $\Psi = \Psi'$. We define 2-*monomorphisms* dually.[2]

Given an enhanced factorization system $(\mathcal{E}, \mathcal{M})$ on a 2-category $A$, we say that $\mathcal{E}$ *consists of* 2-*epimorphisms* if all 1-cells in $\mathcal{E}$ are 2-epimorphisms, and $\mathcal{M}$ *consists of* 2-*monomorphisms* if all of the 1-cells in $\mathcal{M}$ are 2-monomorphisms.

2.7. **Remark.** Another way of characterizing 2-epimorphisms is the following. Consider a 1-cell $\varepsilon \colon F \longrightarrow G$ in a 2-category $A$. To say that $\varepsilon$ is a 2-epimorphism is equivalent to saying that for all objects $H \in A$, the functor between Hom-categories

$$\varepsilon^{\star} = A(\varepsilon, H) \colon A(G, H) \longrightarrow A(F, H)$$

given by pre-composition by $\varepsilon$ is *faithful*. Dually, to say that a 1-cell $\mu \colon G \longrightarrow H$ is a 2-monomorphism is equivalent to saying that for every object $F \in A$, the functor between Hom-categories

$$\mu_{\star} = A(F, \mu) \colon A(F, G) \longrightarrow A(F, H)$$

given by post-composition by $\mu$ is *faithful*. Motivated by this, Dupont and Vitale [2, Def. 3.1] call what we call a 2-monomorphism a "faithful 1-cell" and what we call a 2-epimorphism, a "cofaithful 1-cell". However, Dupont and Vitale's terminology neither seems to be standard nor widely used, and we find our terminology more evocative and clear for our purposes, which is why we have chosen to use it.

Now we are ready to state the main results of this paper.

2.8. **Proposition.** *Suppose that $D$ is a 2-category with an enhanced factorization system $(\mathcal{E}, \mathcal{M})$ and that $C$ is a small 2-category. If $(\mathcal{E}, \mathcal{M})$ separates parallel pairs, $\mathcal{E}$ consists of 2-epimorphisms, and $\mathcal{M}$ consists of 2-monomorphisms, then the levelwise factorization $(\mathcal{E}^C, \mathcal{M}^C)$ produced in Proposition 1.9 defines an enhanced factorization system on $D^C$.*

2.9. **Theorem** (Main Result)**.** *Suppose that $C$ is a small 2-category and that $D$ is a 2-category with an enhanced factorization system $(\mathcal{E}, \mathcal{M})$, and assume the hypotheses of Proposition 2.8. If, in addition, post-composition with 1-cells in $\mathcal{M}$ creates invertible 2-cells, then the levelwise enhanced factorization system $(\mathcal{E}^C, \mathcal{M}^C)$ on $D^C$ of Proposition 2.8 defines an rigid enhanced factorization system on $D^C$.*

Sections 3 and 4 are dedicated to proving Proposition 2.8 and Theorem 2.9. The point of the conditions stated in Proposition 2.8 and Theorem 2.9 is the following: when we were proving Proposition 1.9, we could use the properties of higher cells in $D$ to show that a collection of 1-cells defined levelwise was "coherent" in the sense that the collection defined a 2-natural transformation. However, when we are given a collection of 2-cells in $D$, in general there is not a higher structure to ensure that collection is "coherent" in the sense that the collection defines a modification. Hence, in order to have such a condition hold, we need some extra underlying structure on the enhanced factorization system.

---

[2]This refers to the dual $A^{op}$ given by reversing 1-cells, but not 2-cells.

The following corollary of Theorem 2.9 is a direct application of Lack's result (Theorem 1.8), and the main result of interest to Guillou, May, Merling, and Osorno.

2.9.1. **Corollary.** *Suppose that $C$ is a small 2-category, $D$ is a 2-category with an enhanced factorization system $(\mathcal{E}, \mathcal{M})$, and that $T$ is a 2-monad on $D^C$. Also suppose that*

(2.9.1.a) *$(\mathcal{E}, \mathcal{M})$ separates parallel pairs,*
(2.9.1.b) *$\mathcal{E}$ consists of 2-epimorphisms and $\mathcal{M}$ consists of 2-monomorphisms,*
(2.9.1.c) *and post-composition with 1-cells in $\mathcal{M}$ creates invertible 2-cells.*

*If $T$ preserves $\mathcal{E}^C$, then the inclusion $\mathrm{Alg}_T \hookrightarrow \mathrm{Alg}_T^{ps}$ has a left 2-adjoint, and the components of the unit of the adjunction are equivalences in $\mathrm{Alg}_T^{ps}$.*

As stated at the beginning of the section, Lemmas 2.2 and 2.5 show that these coherence results apply to the case that $D = \mathbf{Cat}$.

2.10. **Example.** For any small 2-category $C$, the levelwise factorization $(\mathcal{B}^C, \mathcal{F}^C)$ defines a rigid enhanced factorization system on $\mathbf{Cat}^C$. Moreover, if $T$ is a 2-monad on $\mathbf{Cat}^C$ which preserves the 2-natural transformations whose components are bijective-on-objects functors, then every pseudo-$T$-algebra is equivalent to a strict $T$-algebra. In particular, this implies Power's result [7, Cor. 3.5].

In [7, p. 170], Power comments that it is straightforward to generalize his coherence result for $\mathbf{Cat}^S$, where $S$ is a small set, to functor 2-categories of the form $\mathbf{Cat}^C$, where $C$ is a small 2-category, which is precisely Example 2.10. However, Power says that "it requires a succession of pasting diagrams, and it is not the case of primary interest; so I omit the proof." In the succeeding sections we will see that this comment is misleading on two accounts. First, there are additional properties about the enhanced factorization system on $\mathbf{Cat}$ that need to be verified in order to make this generalization, and, in particular, they do not simply fall out of the framework that Power provides in [7]. Second, at least in the framework that we have set up, this generalization is not nearly as tedious as Power suggests. Moreover, this result is a crucial piece of Guillou, May, Merling, and Osorno's work [3], and is not nearly as uninteresting as Power suggests.

## 3. Lifting Enhanced Factorization Systems

This section is dedicated to proving Proposition 2.8. The first step in this is to show that, under sufficient hypotheses, the levelwise factorization satisfies the property (1.2.b) of Definition 1.2.

3.1. **Lemma.** *Suppose that $D$ is a 2-category with an enhanced factorization system $(\mathcal{E}, \mathcal{M})$ and that $C$ is a small 2-category. If $(\mathcal{E}, \mathcal{M})$ separates parallel pairs, $\mathcal{E}$ consists of 2-epimorphisms, and $\mathcal{M}$ consists of 2-monomorphisms, then for every diagram*

$$\begin{array}{ccc} F & \xrightarrow{\varepsilon} & F' \\ \alpha \downarrow & \swarrow_{\Psi} & \downarrow \alpha' \\ G & \xrightarrow{\mu} & G' \end{array}$$

*in $D^C$, where $\Psi$ is an invertible 2-cell, $\varepsilon \in \mathcal{E}^C$, and $\mu \in \mathcal{M}^C$, in the 2-category $D^C$ there exists a unique pair $(\delta\colon F' \longrightarrow G, \overline{\Psi}\colon \alpha' \Longrightarrow \mu\delta)$ so that $\delta\varepsilon = \alpha$ and $\overline{\Psi}\varepsilon = \Psi$. Moreover, $\overline{\Psi}$ is necessarily invertible.*

Before proving Lemma 3.1, let us first state a technical sublemma that we use in the proof of Lemma 3.1. The idea of the following result is that if we have a collection of 1-cells $(\beta_c)_{c\in C}$ a 2-category $D$ that satisfy the 1-categorical naturality conditions needed to define a 2-natural transformation between a parallel pair of 2-functors $C \rightrightarrows D$, and, in addition, we know that for each $c \in C$, the 1-cell $\beta_c$ is coherently isomorphic to a 1-cell defining component of an actual 2-natural transformation, then the 1-cells $(\beta_c)_{c\in C}$ themselves satisfy the 2-categorical condition needed to define a 2-natural transformation.

3.1.1. **Sublemma.** *Let $C$ and $D$ be 2-categories, where $C$ is small, let $F, G \colon C \rightrightarrows D$ be a parallel pair of 2-functors, and let $\alpha \colon F \Longrightarrow G$ be a 2-natural transformation. Suppose that we are given*

(3.1.1.a) *a collection of 1-cells $(\beta_c \colon F(c) \longrightarrow G(c))_{c\in C}$ in $D$ satisfying the 1-naturality condition that for each 1-cell $f \colon c \longrightarrow c'$ in $C$, the square*

$$\begin{array}{ccc} F(c) & \xrightarrow{F(f)} & F(c') \\ \beta_c \downarrow & & \downarrow \beta_{c'} \\ G(c) & \xrightarrow[G(f)]{} & G(c') \end{array}$$

*commutes,*

(3.1.1.b) *and a collection of invertible 2-cells $(\Phi_c \colon \alpha_c \Longrightarrow \beta_c)_{c\in C}$ in $D$ satisfying the "modification condition" that for each 1-cell $f \colon c \longrightarrow c'$ of $C$, there is an equality of pasting diagrams*

$$F(c) \underset{\beta_c}{\overset{\alpha_c}{\rightrightarrows}} \Downarrow\Phi_c\ G(c) \xrightarrow{G(f)} G(c') \quad = \quad F(c) \xrightarrow{F(f)} F(c') \underset{\beta_{c'}}{\overset{\alpha_{c'}}{\rightrightarrows}} \Downarrow\Phi_{c'}\ G(c') .$$

*Then the 1-cells $(\beta_c)_{c\in C}$ define the components of a 2-natural transformation $\beta \colon F \Longrightarrow G$. In this case, the invertible 2-cells $(\Phi_c)_{c\in C}$ define the components of an invertible modification from $\alpha$ and $\beta$.*

Since the proof of Sublemma 3.1.1 is a little technical and is a significant detour from the main thrust of the paper, we defer it until Appendix A. We now proceed with the proof of Lemma 3.1, taking Sublemma 3.1.1 for granted.

*Proof of Lemma 3.1.* Since $\varepsilon \in \mathcal{E}^C$, $\mu \in \mathcal{M}^C$, and $(\mathcal{E}, \mathcal{M})$ is an enhanced factorization system on $D$, for each $c \in C$, we have a factorization

$$\begin{array}{ccc} F(c) \xrightarrow{\varepsilon_c} F'(c) & & F(c) \xrightarrow{\varepsilon_c} F'(c) \\ \alpha_c \downarrow\ \ \Psi_c\ \ \downarrow \alpha'_c & = & \alpha_c \downarrow\ \ \overset{\delta_c}{\underset{\widetilde{\Psi}_c}{\Longleftarrow}}\ \ \downarrow \alpha'_c \\ G(c) \xrightarrow[\mu_c]{} G'(c) & & G(c) \xrightarrow[\mu_c]{} G'(c) . \end{array}$$

The goal is to show that the 1-cells $(\delta_c)_{c\in C}$ define the components of a 2-natural transformation $\delta \colon F' \longrightarrow G$, and that the 2-cells $(\widetilde{\Psi}_c)_{c\in C}$ define the components of an invertible modification $\widetilde{\Psi} \colon \alpha' \Longrightarrow \mu\delta$.

*A priori*, the 1-cells $(\delta_c)_{c\in C}$ are not related to one-another. The idea is to show that the 1-cells $(\mu_c\delta_c)_{c\in C}$ and 2-cells $(\widetilde{\Psi}_c)_{c\in C}$ satisfy the conditions of Sublemma 3.1.1, which will prove that the 1-cells $(\delta_c)_{c\in C}$ define the components of a 2-natural transformation, and that the 2-cells $(\widetilde{\Psi}_c)_{c \in C}$ define the components of an invertible modification. In order to apply Sublemma 3.1.1, we need to show that the $(\mu_c \delta_c)_{c \in C}$ satisfy the 1-categorical naturality assumption (3.1.1.a) of Sublemma 3.1.1. Since $\mu$ is a 2-natural transformation, it suffices to show that for each morphism $f : c \longrightarrow c'$ of $C$, the square

$$\begin{array}{ccc} F'(c) & \xrightarrow{F'(f)} & F'(c') \\ \delta_c \downarrow & & \downarrow \delta_{c'} \\ G(c) & \xrightarrow{G(f)} & G(c') \end{array} \qquad (3.1.1)$$

commutes, and then we can apply the 2-naturality of $\mu$ to verify the assumptions of Sublemma 3.1.1. To show that the square (3.1.1) commutes, we exploit the fact that the enhanced factorization system on $D$ separates parallel pairs: it suffices to show that for all morphisms $f : c \longrightarrow c'$ in $C$ we have

$$\delta_{c'} F'(f) \circ \varepsilon_c = G(f) \delta_c \circ \varepsilon_c \quad \text{and} \quad \mu_{c'} \circ \delta_{c'} F'(f) = \mu_{c'} \circ G(f) \delta_c \,.$$

First let us show that $\delta_{c'} F'(f) \circ \varepsilon_c = G(f) \delta_c \circ \varepsilon_c$. Since $\delta_c \varepsilon_c = \alpha_c$ for each $c \in C$, and both $\alpha$ and $\varepsilon$ are 2-natural transformations, for all morphisms $f : c \longrightarrow c'$ in $C$ we have:

$$\begin{array}{ccc} F(c) \xrightarrow{F(f)} F(c') \\ \alpha_c \downarrow \qquad \downarrow \alpha_{c'} \\ G(c) \xrightarrow{G(f)} G(c') \end{array} \;=\; \begin{array}{c} F(c) \xrightarrow{F(f)} F(c') \\ \varepsilon_c \downarrow \qquad \downarrow \varepsilon_{c'} \\ F'(c) \qquad F'(c') \\ \delta_c \downarrow \qquad \downarrow \delta_{c'} \\ G(c) \xrightarrow{G(f)} G(c') \end{array} \;=\; \begin{array}{c} F(c) \\ \varepsilon_c \downarrow \\ F'(c) \xrightarrow{F'(f)} F'(c') \\ \delta_c \downarrow \qquad \downarrow \delta_{c'} \\ G(c) \xrightarrow{G(f)} G(c') \,. \end{array}$$

Hence for all $f : c \longrightarrow c'$ in $C$ we have $\delta_{c'} F'(f) \circ \varepsilon_c = G(f) \delta_c \circ \varepsilon_c$.

Now let us show that for all 1-cells $f : c \longrightarrow c'$ in $C$, the diagram

$$\begin{array}{ccccc} F'(c) & \xrightarrow{F'(f)} & F'(c') & & \\ \delta_c \downarrow & & \downarrow \delta_{c'} & & \\ G(c) & \xrightarrow{G(f)} & G(c') & \xrightarrow{\mu_{c'}} & G'(c') \,. \end{array} \qquad (3.1.2)$$

commutes. To see this, notice that by the 2-naturality of $\mu$ we have $G'(f) \mu_c \delta_c = \mu_{c'} G(f) \delta_c$, so it suffices to show that the invertible 2-cell

$$\widetilde{\Psi}_{c'} F'(f) G'(f) \widetilde{\Psi}_c^{-1} : G'(f) \mu_c \delta_c \Longrightarrow \mu_{c'} \delta_{c'} F'(f)$$

is an identity. This follows from the fact that $\mathcal{E}$ consists of 2-epimorphisms: whiskering with $\varepsilon_c$ and applying the fact that $\Psi^{-1}$ is a modification we see that

$$G'(f) \widetilde{\Psi}_c^{-1} \varepsilon_c = G'(f) \Psi_c^{-1} = \Psi_{c'}^{-1} F(f) \,.$$

By the 2-naturality of $\varepsilon$, we have

$$\widetilde{\Psi}_{c'}^{-1} F'(f) \varepsilon_c = \widetilde{\Psi}_{c'}^{-1} \varepsilon_{c'} F(f) = \Psi_{c'}^{-1} F(f) \,.$$

Finally, since $\varepsilon_c$ is a 2-epimorphism, we see that

$$(\widetilde{\Psi}_{c'} F'(f))^{-1} = \widetilde{\Psi}_{c'}^{-1} F'(f) = G'(f) \widetilde{\Psi}_c^{-1} \,, \qquad (3.1.3)$$

so $\widetilde{\Psi}_{c'} F'(f) G'(f) \widetilde{\Psi}_c^{-1}$ is an identity, as desired. Then since $\delta_{c'} F'(f) \circ \varepsilon_c = G(f) \delta_c \circ \varepsilon_c$ and the enhanced factorization system on $D$ separates parallel pairs, the square (3.1.1) commutes.

The commutativity of the square (3.1.1) and the 2-naturality of $\mu$ together imply that for every 1-cell $f\colon c \longrightarrow c'$ in $\mathbf{C}$, the square

$$\begin{array}{ccc} F'(c) & \xrightarrow{F'(f)} & F'(c') \\ {\scriptstyle \mu_c \delta_c}\downarrow & & \downarrow{\scriptstyle \mu_{c'} \delta_{c'}} \\ G(c) & \xrightarrow[G(f)]{} & G(c') \end{array} \qquad (3.1.4)$$

commutes.

Now let us show that the invertible 2-cells $\widetilde{\Psi}_c$ in $\mathbf{D}$ satisfy the "modification condition" that for each 1-cell $f\colon c \longrightarrow c'$ of $\mathbf{C}$, there is an equality of pasting diagrams

$$F'(c) \; \overset{\alpha'_c}{\underset{\mu_c \delta_c}{\widetilde{\Psi}_c\Downarrow}} \; G'(c) \xrightarrow{G'(f)} G'(c') \quad = \quad F'(c) \xrightarrow{F'(f)} F'(c') \; \overset{\alpha'_{c'}}{\underset{\mu_{c'} \delta_{c'}}{\Downarrow\widetilde{\Psi}_{c'}}} \; G'(c') \, . \qquad (3.1.5)$$

However, this follows immediately from equation (3.1.3), using the fact that the modifications appearing in equation (3.1.3) are invertible and that the inverse of $G'(f)\widetilde{\Psi}_c^{-1}$ is given by $G'(f)\widetilde{\Psi}_c$.

The commutativity of the square (3.1.4) and the equality of pasting diagrams (3.1.5) for every 1-cell $f\colon c \longrightarrow c'$ in $\mathbf{C}$ puts us in the situation where we can apply Sublemma 3.1.1. Sublemma 3.1.1 shows that whenever we are given the data of objects, 1-cells, and 2-cells in $\mathbf{C}$, as displayed below:

$$c \; \overset{f}{\underset{g}{\Downarrow\lambda}} \; c' \, ,$$

we have an equality of pasting diagrams

$$F'(c) \; \overset{F'(f)}{\underset{F'(g)}{F'(\lambda)\Downarrow}} \; F'(c') \xrightarrow{\mu_{c'} \delta_{c'}} G'(c') \quad = \quad F'(c) \xrightarrow{\mu_c \delta_c} G'(c) \; \overset{G'(f)}{\underset{G'(g)}{\Downarrow G'(\lambda)}} \; G'(c') \, .$$

By the 2-naturality of $\mu$, we see that

$$F'(c) \; \overset{F'(f)}{\underset{F'(g)}{F'(\lambda)\Downarrow}} \; F'(c') \xrightarrow{\mu_{c'} \delta_{c'}} G'(c') \quad = \quad F'(c) \xrightarrow{\delta_c} G(c) \; \overset{G(f)}{\underset{G(g)}{\Downarrow G(\lambda)}} \; G(c') \xrightarrow{\mu_{c'}} G'(c') \, .$$

Since $\mu_{c'}$ is a 2-monomorphism, this implies that

$$F'(c) \; \overset{F'(f)}{\underset{F'(g)}{F'(\lambda)\Downarrow}} \; F'(c') \xrightarrow{\delta_{c'}} G(c') \quad = \quad F'(c) \xrightarrow{\delta_c} G(c) \; \overset{G(f)}{\underset{G(g)}{\Downarrow G(\lambda)}} \; G(c') \, ,$$

which shows that the 1-cells $(\delta_c)_{c \in \mathbf{C}}$ satisfy the 2-categorical condition to define a 2-natural transformation, hence define the components of a 2-natural transformation $\delta$. The uniqueness of $\delta$ comes for free from the uniqueness property of the components $(\delta_c)_{c \in \mathbf{C}}$. □

3.2. **Remark.** By analyzing the proof of Lemma 3.1, it is easy to see that in the case that $C$ is an ordinary 1-category, the assumption that $\mathcal{M}$ consists of 2-monomorphisms is superfluous. However, to handle the rigidity hypothesis of Theorem 1.8, in §4 we assume that $\mathcal{M}$ consists of 2-monomorphisms. In the specialized case that $C$ is a 1-category, to apply our results from §4 the assumption that $\mathcal{M}$ consists of 2-monomorphisms is still necessary, hence, if coherence results is what one is after, then there is no harm assuming this from the onset.

3.3. **Lemma.** *Suppose that $D$ is a 2-category with an enhanced factorization system $(\mathcal{E}, \mathcal{M})$ and that $C$ is a small 2-category. If $(\mathcal{E}, \mathcal{M})$ separates parallel pairs, $\mathcal{E}$ consists of 2-epimorphisms, and $\mathcal{M}$ consists of 2-monomorphisms then the levelwise factorization $(\mathcal{E}^C, \mathcal{M}^C)$ satisfies the last condition (1.2.c) of Definition 1.2, hence defines an enhanced factorization system on $D^C$.*

*Proof.* Suppose that we are in the situation indicated in the last condition (1.2.c) of Definition 1.2 (since it is lengthy, we will not spell it out again here.) For each $c \in C$, write $\delta_{1,c}$ and $\delta_{2,c}$ for the components of $\delta_1$ at $c$ and $\delta_2$ at $c$, respectively. Since $(\mathcal{E}, \mathcal{M})$ is an enhanced factorization system on $D$, for each $c \in C$, there exists a unique 2-cell $\Delta_c \colon \delta_{1,c} \Longrightarrow \delta_{2,c}$ in $D$ so that $\Delta_c \varepsilon_c = \Phi_c$ and $\mu_c \Delta_c = \Phi'_c$. To see that the 2-cells $(\Delta_c)_{c \in C}$ assemble into an modification $\Delta \colon \delta_1 \Longrightarrow \delta_2$, notice that for all morphisms $f \colon c \longrightarrow c'$ in $C$, since $\Phi$ is a modification we have

$$G(f)\Delta_c \varepsilon_c = G(f)\Phi_c = \Phi_{c'} F(f) \,.$$

Similarly, since $\varepsilon$ is a 2-natural transformation and $\Delta_{c'}\varepsilon_{c'} = \Phi_{c'}$, we see that

$$\Delta_{c'} F'(f)\varepsilon_c = \Delta_{c'}\varepsilon_{c'} F(f) = \Phi_{c'} F(f) \,,$$

hence $G(f)\Delta_c \varepsilon_c = \Delta_{c'} F'(f)\varepsilon_c$. Then since $\mathcal{E}$ consists of 2-epimorphisms, we also see that $G(f)\Delta_c = \Delta_{c'} F'(f)$. Hence the 2-cells $(\Delta_c)_{c \in C}$ satisfy the necessary conditions to define a modification. □

Combining Proposition 1.9 and Lemmas 3.1 and 3.3 proves Proposition 2.8.

## 4. Conditions for Rigidity

In this section we analyze conditions on the enhanced factorization system on $D$ which yield a rigid enhanced factorization system on $D^C$, for any small 2-category $C$. Again, since the rigidity condition does not have a uniqueness statement, we have no way of guaranteeing 2-cells in the rigidity condition assemble coherently to define a modification, hence we must impose an additional condition to guarantee this coherence. This can be achieved by assuring that the 2-cells specified by the rigidity of the enhanced factorization system on $D$ can be chosen "canonically"; assuming that post-composition with 1-cells in $\mathcal{M}$ creates invertible 2-cells allows us to do this.

4.1. **Lemma.** *Suppose that $A$ is a 2-category with an enhanced factorization system $(\mathcal{E}, \mathcal{M})$. If post-composition with 1-cells in $\mathcal{M}$ creates invertible 2-cells, then given 1-cells $\mu \colon F \rightleftarrows G : \alpha$, where $\mu \in \mathcal{M}$, and an invertible 2-cell $\Psi \colon \mu\alpha \Longrightarrow \mathrm{id}_G$, there exists a unique invertible 2-cell $\widehat{\Phi} \colon \alpha\mu \Longrightarrow \mathrm{id}_F$ so that $\mu\widehat{\Phi} = \Psi\mu$. In particular, the enhanced factorization system on $A$ is rigid.*

*Proof.* Suppose that we have 1-cells $\mu \colon F \rightleftarrows G : \alpha$, where $\mu \in \mathcal{M}$, and we have an invertible 2-cell $\Psi \colon \mu\alpha \Longrightarrow \mathrm{id}_G$. Then the 2-cell $\Psi\mu \colon \mu\alpha\mu \Longrightarrow \mu$ is invertible. For notational simplicity, write $\Phi := \Psi\mu$. Then since post-composition with $\mu$ creates invertible 2-cells, there is a unique invertible 2-cell

$$\widehat{\Phi} \colon \alpha\mu \Longrightarrow \mathrm{id}_F$$

so that $\mu\widehat{\Phi} = \Psi\mu$. In particular, $\alpha\mu \cong \text{id}_F$. □

**4.2. Lemma.** *Suppose that $\mathbf{D}$ is a 2-category with an enhanced factorization system $(\mathcal{E}, \mathcal{M})$ and that $\mathbf{C}$ is a small 2-category. If post-composition with 1-cells in $\mathcal{M}$ creates invertible 2-cells and $\mathcal{M}$ consists of 2-monomorphsims, then post-composition with 1-cells in $\mathcal{M}^{\mathbf{C}}$ creates invertible 2-cells in $\mathbf{D}^{\mathbf{C}}$.*

*Proof.* Suppose that we are given 1-cells $\mu\colon F \rightleftarrows G : \alpha$, where $\mu \in \mathcal{M}^{\mathbf{C}}$, and we have an invertible 2-cell $\Psi\colon \mu\alpha \implies \text{id}_G$. Then for each $c \in \mathbf{C}$, there exists an invertible 2-cell $\Psi_c \colon \mu_c\alpha_c \implies \text{id}_{G(c)}$ in $\mathbf{D}$. By Lemma 4.1, there exists a unique invertible 2-cell

$$\widehat{\Phi}_c \colon \alpha_c\mu_c \implies \text{id}_{F(c)}$$

so that $\mu_c\widehat{\Phi}_c = \Psi_c\mu_c$. We want to show that the $\widehat{\Phi}_c$ assemble into the components of an invertible modification.

To do this, suppose that $f\colon c \longrightarrow c'$ is a morphism in $C$. Writing $\Phi \coloneqq \Psi\mu$, since $\Phi$ is a modification, we have

$$F(c) \xrightarrow[\mu_c]{\overset{\mu_c\alpha_c\mu_c}{\Phi_c\Downarrow}} G(c) \xrightarrow{G(f)} G(c') \quad = \quad F(c) \xrightarrow{F(f)} F(c') \xrightarrow[\mu_{c'}]{\overset{\mu_{c'}\alpha_{c'}\mu_{c'}}{\Downarrow\Phi_{c'}}} G(c')\;.$$

Hence we see that

$$\mu_{c'}\widehat{\Phi}_{c'}F(f) = \Phi_{c'}F(f) = G(f)\Phi_c\;.$$

By the naturality of $\mu$ and the fact that $\mu\widehat{\Phi} = \Phi$ we see that

$$\mu_{c'}F(f)\widehat{\Phi}_c = G(f)\mu_c\widehat{\Phi}_c = G(f)\Phi_c\;.$$

Then since

$$\mu_{c'}\widehat{\Phi}_{c'}F(f) = \mu_{c'}F(f)\widehat{\Phi}_c$$

and $\mu_{c'}$ is a 2-monomorphsim, we have $\widehat{\Phi}_{c'}F(f) = F(f)\widehat{\Phi}_c$, hence the 2-cells $(\widehat{\Phi}_c)_{c\in \mathbf{C}}$ in $\mathbf{D}$ assemble into an invertible modification $\widehat{\Phi}$. □

**4.2.1. Corollary.** *Suppose that $\mathbf{D}$ is a 2-category with an enhanced factorization system $(\mathcal{E}, \mathcal{M})$ and that $\mathbf{C}$ is a small 2-category. If*

(4.2.1.a) *$(\mathcal{E}, \mathcal{M})$ separates parallel pairs,*
(4.2.1.b) *$\mathcal{E}$ consists of 2-epimorphisms and $\mathcal{M}$ consists of 2-monomorphisms,*
(4.2.1.c) *and post-composition with 1-cells in $\mathcal{M}$ creates invertible 2-cells,*

*then the enhanced factorization system $(\mathcal{E}^{\mathbf{C}}, \mathcal{M}^{\mathbf{C}})$ on $\mathbf{D}^{\mathbf{C}}$ is rigid.*

Hence Proposition 2.8 and corollary 4.2.1 together prove the main result, Theorem 2.9.

## Appendix A. The Proof of a Technical Result

In this appendix we recall and prove Sublemma 3.1.1.

**3.1.1. Sublemma.** *Let $\mathbf{C}$ and $\mathbf{D}$ be 2-categories, where $\mathbf{C}$ is small, let $F, G\colon \mathbf{C} \rightrightarrows \mathbf{D}$ be a parallel pair of 2-functors, and let $\alpha\colon F \implies G$ be a 2-natural transformation. Suppose that we are given* (3.1.1.a) *a collection of 1-cells $(\beta_c \colon F(c) \longrightarrow G(c))_{c \in \mathbf{C}}$ in $\mathbf{D}$ satisfying the 1-naturality condition that for each 1-cell $f \colon c \longrightarrow c'$ in $\mathbf{C}$, the square*

$$\begin{array}{ccc} F(c) & \xrightarrow{F(f)} & F(c') \\ \beta_c \downarrow & & \downarrow \beta_{c'} \\ G(c) & \xrightarrow{G(f)} & G(c') \end{array}$$

*commutes,*

(3.1.1.b) *and a collection of invertible 2-cells $(\Phi_c \colon \alpha_c \Longrightarrow \beta_c)_{c \in \mathbf{C}}$ in $\mathbf{D}$ satisfying the "modification condition" that for each 1-cell $f \colon c \longrightarrow c'$ of $\mathbf{C}$, there is an equality of pasting diagrams*

$$F(c) \overset{\alpha_c}{\underset{\beta_c}{\overset{\Phi_c \Downarrow}{\rightrightarrows}}} G(c) \xrightarrow{G(f)} G(c') \quad = \quad F(c) \xrightarrow{F(f)} F(c') \overset{\alpha_{c'}}{\underset{\beta_{c'}}{\overset{\Downarrow \Phi_{c'}}{\rightrightarrows}}} G(c').$$

*Then the 1-cells $(\beta_c)_{c \in \mathbf{C}}$ define the components of a 2-natural transformation $\beta \colon F \Longrightarrow G$. In this case, the invertible 2-cells $(\Phi_c)_{c \in \mathbf{C}}$ define the components of an invertible modification from $\alpha$ and $\beta$.*

*Proof of Sublemma 3.1.1.* In order to prove that the collection of 1-cells $(\beta_c)_{c \in \mathbf{C}}$ defines a 2-natural transformation, we need to show that whenever we are given the data of of objects, 1-cells, and a 2-cell in $\mathbf{C}$, as displayed below

$$c \overset{f}{\underset{g}{\overset{\Downarrow \lambda}{\rightrightarrows}}} c',$$

we have an equality of pasting diagrams

$$F(c) \overset{F(f)}{\underset{F(g)}{\overset{F(\lambda) \Downarrow}{\rightrightarrows}}} F(c') \xrightarrow{\beta_{c'}} G(c') \quad = \quad F(c) \xrightarrow{\beta_c} G(c) \overset{G(f)}{\underset{G(g)}{\overset{\Downarrow G(\lambda)}{\rightrightarrows}}} G(c'). \qquad (\text{A.1.1})$$

Since $\alpha$ is a 2-natural transformation we know that

$$F(c) \overset{F(f)}{\underset{F(g)}{\overset{F(\lambda) \Downarrow}{\rightrightarrows}}} F(c') \xrightarrow{\alpha_{c'}} G(c') \quad = \quad F(c) \xrightarrow{\alpha_c} G(c) \overset{G(f)}{\underset{G(g)}{\overset{\Downarrow G(\lambda)}{\rightrightarrows}}} G(c'),$$

so the idea is to pre-compose the left-hand diagram in (A.1.1) with the invertible 2-cell $\Phi_{c'} F(f)$, and show that this composite is equal to the right-hand diagram in (A.1.1) pre-composed with $\Phi_{c'} F(f)$, proving the equality indicated in (A.1.1).

First notice that by the uniqueness of the definition of horizontal composite $\Phi_{c'} \star F(\lambda)$ displayed below

$$F(c) \overset{F(f)}{\underset{F(g)}{\overset{F(\lambda) \Downarrow}{\rightrightarrows}}} F(c') \overset{\alpha_{c'}}{\underset{\beta_{c'}}{\overset{\Downarrow \Phi_{c'}}{\rightrightarrows}}} G(c'),$$

$\Phi_{c'} \star F(\lambda)$ is equal to both of the vertical composites $\beta_{c'} F(\lambda) \circ \Phi_{c'} F(f)$ and $\Phi_{c'} F(g) \circ \alpha_{c'} F(\lambda)$. Using this fact, the fact that

$$F(c) \xrightarrow[\beta_c]{\alpha_c} \Downarrow\Phi_c \; G(c) \xrightarrow{G(g)} G(c') \quad = \quad F(c) \xrightarrow{F(g)} F(c') \xrightarrow[\beta_{c'}]{\alpha_{c'}} \Downarrow\Phi_{c'} \; G(c'),$$

and the fact that the 1-cells $(\beta_c)_{c \in C}$ satisfy the 1-naturality condition (3.1.1.a) of Sublemma 3.1.1, by the 2-naturality of $\alpha$ we see that

$$F(c) \xrightarrow[F(g)]{F(f)} \Downarrow F(\lambda) \; F(c') \xrightarrow[\beta_{c'}]{\alpha_{c'}} \Downarrow\Phi_{c'} \; G(c') \quad = \quad F(c) \xrightarrow[\beta_c]{\alpha_c} \Downarrow\Phi_c \; G(c) \xrightarrow[G(g)]{G(f)} \Downarrow G(\lambda) \; G(c').$$

Then by the fact that

$$F(c) \xrightarrow[\beta_c]{\alpha_c} \Downarrow\Phi_c \; G(c) \xrightarrow{G(f)} G(c') \quad = \quad F(c) \xrightarrow{F(f)} F(c') \xrightarrow[\beta_{c'}]{\alpha_{c'}} \Downarrow\Phi_{c'} \; G(c'),$$

and the definition of the horizontal composite, we see that

$$F(c) \xrightarrow[\beta_c]{\alpha_c} \Downarrow\Phi_c \; G(c) \xrightarrow[G(g)]{G(f)} \Downarrow G(\lambda) \; G(c') \quad = \quad F(c) \xrightarrow[\beta_c]{\alpha_{c'} F(f)}_{\Phi_{c'} F(f) \Downarrow} G(c) \xrightarrow[G(g)]{G(f)} \Downarrow G(\lambda) \; G(c').$$

This proves that

$$F(c) \xrightarrow[F(g)]{F(f)} \Downarrow F(\lambda) \; F(c') \xrightarrow[\beta_{c'}]{\alpha_{c'}} \Downarrow\Phi_{c'} \; G(c') \quad = \quad F(c) \xrightarrow[\beta_c]{\alpha_{c'} F(f)}_{\Phi_{c'} F(f) \Downarrow} G(c) \xrightarrow[G(g)]{G(f)} \Downarrow G(\lambda) \; G(c'),$$

and since the 2-cell $\Phi_{c'} F(f)$ is invertible, this proves that the equality of pasting diagrams displayed in (A.1.1) holds. Thus the 1-cells $(\beta_c)_{c \in C}$ define the components of a 2-natural transformation, as claimed. Finally, the fact that the invertible 2-cells $(\Phi_c)_{c \in C}$ define the components of an invertible modification from $\alpha$ and $\beta$ follows immediately from condition (3.1.1.b). □

Massachusetts Institute of Technology, Department of Mathematics, 77 Massachusetts Avenue, Cambridge, MA 02139-4307, USA
  *E-mail address*: phaine@mit.edu